# RANK-BASED ESTIMATION FOR ALL-PASS TIME SERIES MODELS

By Beth Andrews,[1] Richard A. Davis[1,2] and F. Jay Breidt[2]

*Northwestern University, Colorado State University
and Colorado State University*

An autoregressive-moving average model in which all roots of the autoregressive polynomial are reciprocals of roots of the moving average polynomial and vice versa is called an *all-pass time series model*. All-pass models are useful for identifying and modeling noncausal and noninvertible autoregressive-moving average processes. We establish asymptotic normality and consistency for rank-based estimators of all-pass model parameters. The estimators are obtained by minimizing the rank-based residual dispersion function given by Jaeckel [*Ann. Math. Statist.* **43** (1972) 1449–1458]. These estimators can have the same asymptotic efficiency as maximum likelihood estimators and are robust. The behavior of the estimators for finite samples is studied via simulation and rank estimation is used in the deconvolution of a simulated water gun seismogram.

**1. Introduction.** Autoregressive-moving average (ARMA) models, the standard linear time series models for stationary data, are often fit to observed series using Gaussian likelihood, least-squares, or related second-order moment estimation techniques. These are effective methods for finding fitted ARMA models with second-order moment properties that resemble those of an observed series, whether or not the data are Gaussian. However, because every Gaussian ARMA process has a causal, invertible ARMA representation (all roots of the autoregressive and moving average polynomials are outside the unit circle), in the non-Gaussian case, the second-order methods are unable to identify a noncausal (at least one root of the autoregressive polynomial is inside the unit circle) or noninvertible (at least one root of

Received January 2005; revised July 2006.

[1]Supported in part by NSF Grants DMS-99-72015 and DMS-03-08109.

[2]Supported in part by EPA STAR Grant CR-829095.

the moving average polynomial is inside the unit circle) ARMA series. Fitted ARMA models obtained using second-order techniques may not, therefore, most effectively capture the higher-order moment structure of the data. Consequently, an effort to identify noncausal and noninvertible series should be part of any ARMA fitting procedure. In this paper, we discuss all-pass models, which are useful tools for identifying and modeling noncausal and noninvertible ARMA processes.

All-pass models are ARMA models in which the roots of the autoregressive polynomial are reciprocals of roots of the moving average polynomial and vice versa. These models generate uncorrelated (white noise) time series that are not independent in the non-Gaussian case. As discussed in [2], an all-pass series can be obtained by fitting a causal, invertible ARMA model to a series generated by a causal, noninvertible ARMA model. The residuals follow an all-pass model of order $r$, where $r$ is the number of roots of the true moving average polynomial inside the unit circle. Consequently, by identifying the all-pass order of the residuals, the order of noninvertibility of the ARMA model can be determined without considering all possible configurations of roots inside and outside the unit circle, which is computationally prohibitive for large-order models. Noninvertible ARMA models have appeared, for example, in vocal tract filters [8, 9], in the analysis of unemployment rates [13] and in seismogram deconvolution [2, 19]. All-pass models can be used similarly to fit noncausal ARMA models [6]. See [6] for a list of applications for noncausal models.

Estimation methods based on second-order moment techniques cannot identify all-pass models because Gaussian all-pass series are independent. Thus, cumulant-based estimators, using cumulants of order greater than two, are often used to estimate these models [8, 9, 11]. Breidt, Davis and Trindade [6] consider a least absolute deviations (LAD) estimation approach which is motivated by the likelihood of an all-pass model with Laplace (two-sided exponential) noise, and Andrews, Davis and Breidt [2] consider a maximum likelihood (ML) estimation approach. The LAD and ML estimators are consistent and asymptotically normal. However, the LAD estimation procedure is limited by the assumption that the mean and median for the noise are equivalent, and the ML procedure is limited by the assumption that the probability density function for the noise is symmetric and known to within some parameter values.

In this paper, we consider a rank-based estimation technique first proposed by Jaeckel [14] for estimating linear regression parameters. Jaeckel's estimator minimizes the sum of model residuals weighted by a function of residual rank. We study the asymptotic properties of Jaeckel's rank ($R$) estimator in the case of all-pass parameter estimation. This $R$-estimator is more robust than the LAD and ML estimators; it is consistent and asymptotically normal under less stringent conditions. In addition, when $R$-estimation is used in lieu of LAD or ML, efficiency need not be sacrificed. There exists a weight function for which $R$-estimation is asymptotically equivalent to LAD estimation and, when the noise distribution is known, the weight function can be chosen so that $R$-estimation is asymptotically equivalent to ML estimation. We also find that when the Wilcoxon weight function (a linear weight function) is used, $R$-estimation is (relatively) very efficient for a large class of noise distributions. Another advantage of $R$-estimation is that one has the flexibility to choose a weight function that tends to produce relatively smooth $R$-objective functions which can be minimized fairly easily.

Because the objective function for Jaeckel's $R$-estimation method involves not only the residual ranks, but also the residual values, this is not pure $R$-estimation. Koul and Ossiander [16], Koul and Saleh [17], Mukherjee and Bai [20] and Terpstra, McKean and Naranjo [23] consider related rank-based estimation approaches for autoregressive model parameters. Also, Allal, Kaaouachi and Paindaveine [1] examine a pure $R$-estimator for ARMA model parameters based on correlations of weighted residual ranks. The results for this pure $R$-estimator are not applicable to all-pass model parameters because the parameters in the autoregressive polynomial of an all-pass model are functions of parameters in the moving average polynomial and vice versa.

In Section 2 we consider Jaeckel's $R$-function in the context of all-pass parameter estimation. Asymptotic normality for $R$-estimators is established under mild conditions and order selection is discussed in Section 3. Proofs of the lemmas used to establish the results of Section 3 can be found in the Appendix. We study the behavior of the estimators for finite samples via simulation in Section 4.1 and use $R$-estimation in the deconvolution of a simulated water gun seismogram in Section 4.2.

## 2. Preliminaries.

2.1. *All-pass models.* Let $B$ denote the backshift operator ($B^k X_t = X_{t-k}$, $k = 0, \pm 1, \pm 2, \ldots$) and let $\phi(z) = 1 - \phi_1 z - \cdots - \phi_p z^p$ be a $p$th order autoregressive polynomial, where $\phi(z) \neq 0$ for $|z| = 1$. The filter $\phi(B)$ is said to be *causal* if all the roots of $\phi(z)$ are outside the unit circle in the complex plane. In this case, for a sequence $\{W_t\}$, $\phi^{-1}(B)W_t = (\sum_{j=0}^{\infty} \psi_j B^j)W_t = \sum_{j=0}^{\infty} \psi_j W_{t-j}$, a function of only the past and present $\{W_t\}$. If $\phi(B)$ is causal, then the filter $B^p \phi(B^{-1})$ is *purely noncausal* and hence $B^{-p} \phi^{-1}(B^{-1})W_t = (\sum_{j=0}^{\infty} \psi_j B^{-p-j})W_t = \sum_{j=0}^{\infty} \psi_j W_{t+p+j}$, a function of only the present and future $\{W_t\}$. See, for example, Chapter 3 of [7].

Let $\phi_0(z) = 1 - \phi_{01} z - \cdots - \phi_{0p} z^p$, where $\phi_0(z) \neq 0$ for $|z| \leq 1$. Define $\phi_{00} = 1$ and $r = \max\{0 \leq j \leq p : \phi_{0j} \neq 0\}$. Then a causal all-pass time series is the ARMA series $\{X_t\}$ which satisfies the difference equations

$$(2.1) \qquad \phi_0(B)X_t = \frac{B^r \phi_0(B^{-1})}{-\phi_{0r}} Z_t^*$$

or

$$X_t - \phi_{01}X_{t-1} - \cdots - \phi_{0r}X_{t-r} = Z_t^* + \frac{\phi_{0,r-1}}{\phi_{0r}} Z_{t-1}^* + \cdots + \frac{\phi_{01}}{\phi_{0r}} Z_{t-r+1}^* - \frac{1}{\phi_{0r}} Z_{t-r}^*,$$

where the series $\{Z_t^*\}$ is an independent and identically distributed (i.i.d.) sequence of random variables with mean 0, variance $\sigma^2 \in (0, \infty)$ and distribution function $F$. The true order of the all-pass model is $r$ $(0 \le r \le p)$. Observe that the roots of the autoregressive polynomial $\phi_0(z)$ are reciprocals of the roots of the moving average polynomial $-\phi_{0r}^{-1} z^r \phi_0(z^{-1})$ and vice versa.

The spectral density for $\{X_t\}$ in (2.1) is

$$\frac{|e^{-ir\omega}|^2 |\phi_0(e^{i\omega})|^2}{\phi_{0r}^2 |\phi_0(e^{-i\omega})|^2} \frac{\sigma^2}{2\pi} = \frac{\sigma^2}{\phi_{0r}^2 2\pi},$$

which is constant for $\omega \in [-\pi, \pi]$. Thus, $\{X_t\}$ is an uncorrelated sequence. In the case of Gaussian $\{Z_t^*\}$, this implies that $\{X_t\}$ is i.i.d. $\mathrm{N}(0, \sigma^2 \phi_{0r}^{-2})$, but independence does not hold in the non-Gaussian case if $r \ge 1$ (see [5]). The model (2.1) is called all-pass because the power transfer function of the all-pass filter passes all the power for every frequency in the spectrum. In other words, an all-pass filter does not change the distribution of power over the spectrum.

We can express (2.1) as

$$(2.2) \qquad \phi_0(B)X_t = \frac{B^p \phi_0(B^{-1})}{-\phi_{0r}} Z_t,$$

where $\{Z_t\} = \{Z_{t+p-r}^*\}$ is an i.i.d. sequence of random variables with mean 0, variance $\sigma^2$ and distribution function $F$. Rearranging (2.2) and setting $z_t = \phi_{0r}^{-1} Z_t$, we have the backward recursion $z_{t-p} = \phi_{01}z_{t-p+1} + \cdots + \phi_{0p}z_t - (X_t - \phi_{01}X_{t-1} - \cdots - \phi_{0p}X_{t-p})$. An analogous recursion for an arbitrary, causal autoregressive polynomial $\phi(z) = 1 - \phi_1 z - \cdots - \phi_p z^p$ can be defined as

$$(2.3) \qquad z_{t-p}(\boldsymbol{\phi}) = \begin{cases} 0, & t = n+p, \ldots, n+1, \\ \phi_1 z_{t-p+1}(\boldsymbol{\phi}) + \cdots + \phi_p z_t(\boldsymbol{\phi}) - \phi(B)X_t, \\ & t = n, \ldots, p+1, \end{cases}$$

where $\boldsymbol{\phi} := (\phi_1, \ldots, \phi_p)'$. Let $\boldsymbol{\phi}_0 = (\phi_{01}, \ldots, \phi_{0p})' = (\phi_{01}, \ldots, \phi_{0r}, 0, \ldots, 0)'$ denote the true parameter vector and note that $\{z_t(\boldsymbol{\phi}_0)\}_{t=1}^{n-p}$ closely approximates $\{z_t\}_{t=1}^{n-p}$; the error is due to the initialization with zeros. Although $\{z_t\}$ is i.i.d., $\{z_t(\boldsymbol{\phi}_0)\}_{t=1}^{n-p}$ is not i.i.d. if $r \ge 1$.

2.2. *Jaeckel's rank function.* Suppose we have a realization of length $n$, $\{X_t\}_{t=1}^n$, from (2.1). Let $\lambda$ be a function from $(0,1)$ to $\mathbb{R}$ such that

A1. $\lambda$ is strictly increasing and $\lambda(s) = -\lambda(1-s)$ for all $s \in (0,1)$.

If $\phi$ forms a causal $p$th order autoregressive polynomial and $\{R_t(\phi)\}_{t=1}^{n-p}$ contains the ranks of $\{z_t(\phi)\}_{t=1}^{n-p}$ from (2.3), then the $R$-function evaluated at $\phi$ with weight function $\lambda$ is

$$(2.4) \qquad D(\phi) := \sum_{t=1}^{n-p} \lambda\left(\frac{R_t(\phi)}{n-p+1}\right) z_t(\phi).$$

Because it tends to be near zero when the elements of $\{z_t(\phi)\}$ are similar, (2.4) is a measure of the dispersion of the residuals $\{z_t(\phi)\}$. When $\{z_{(t)}(\phi)\}_{t=1}^{n-p}$ is the series $\{z_t(\phi)\}_{t=1}^{n-p}$ ordered from smallest to largest, (2.4) can also be written as $D(\phi) = \sum_{t=1}^{n-p} \lambda(t/(n-p+1)) z_{(t)}(\phi)$. A popular choice for the weight function is $\lambda(s) = s - 1/2$. In this case, the weights $\{\lambda(t/(n-p+1))\}_{t=1}^{n-p}$ are known as *Wilcoxon scores*.

We give some properties for $D$ in the following theorem. Jaeckel [14] shows that the same properties hold for the $R$-function in the linear regression case.

THEOREM 2.1. *Assume* A1 *holds. For any* $\phi \in \mathbb{R}^p$, *if*

$$\{P_1(\phi), \ldots, P_{(n-p)!}(\phi)\}$$
$$= \{\{z_{1,1}(\phi), \ldots, z_{1,n-p}(\phi)\}, \ldots, \{z_{(n-p)!,1}(\phi), \ldots, z_{(n-p)!,n-p}(\phi)\}\}$$

*contains the* $(n-p)!$ *permutations of the sequence* $\{z_t(\phi)\}_{t=1}^{n-p}$, *then*

$$D(\phi) = \sup_{j \in \{1,\ldots,(n-p)!\}} \sum_{t=1}^{n-p} \lambda\left(\frac{t}{n-p+1}\right) z_{j,t}(\phi).$$

*In addition,* $D$ *is a nonnegative, continuous function on* $\mathbb{R}^p$, *and* $D(\phi) = 0$ *if and only if the elements of* $\{z_t(\phi)\}_{t=1}^{n-p}$ *are all equal.*

PROOF. See the proof of Theorem 1 in [14]. □

## 3. Asymptotic results.

3.1. *Parameter estimation.* In order to establish asymptotic normality for $R$-estimators of $\phi_0$, we make the following additional assumptions:

A2. $F$, the distribution function for the noise, is strictly increasing and differentiable on $\mathbb{R}$ with density $f$;

A3. $f$ is uniformly continuous on $\mathbb{R}$ with $\sup_{s \in \mathbb{R}} |s| f(s) < \infty$;

A4. the derivative of the weight function $\lambda$ exists and is uniformly continuous on $(0, 1)$.

Also, let $\tilde{J} = \int_0^1 \lambda^2(s)\,ds$, $\tilde{K} = \int_0^1 F^{-1}(s)\lambda(s)\,ds$ and $\tilde{L} = \int_0^1 f(F^{-1}(s)) \times \lambda'(s)\,ds$ and assume

A5. $\sigma^2 \tilde{L} > \tilde{K}$.

THEOREM 3.1.  *If* A1–A5 *hold, then there exists a sequence of minimizers $\hat{\phi}_R$ of $D(\cdot)$ in (2.4) such that*

$$n^{1/2}(\hat{\phi}_R - \phi_0) \xrightarrow{d} \mathbf{Y} \sim \mathrm{N}(\mathbf{0}, \mathbf{\Sigma}), \tag{3.1}$$

*where* $\mathbf{\Sigma} := (\sigma^2 \tilde{J} - \tilde{K}^2)/[2(\sigma^2 \tilde{L} - \tilde{K})^2]\sigma^2 \mathbf{\Gamma}_p^{-1}$, $\mathbf{\Gamma}_p := [\gamma(j-k)]_{j,k=1}^p$ *and* $\gamma(\cdot)$ *is the autocovariance function for the autoregressive process* $\{(1/\phi_0(B))Z_t\}$.

PROOF.  $D(\phi) - D(\phi_0) = S_n(\sqrt{n}(\phi - \phi_0))$, where $S_n(\cdot)$ is defined in Lemma A.5 of the Appendix. Because $\mathbf{Y} := -|\phi_{0r}|\sigma^2 \mathbf{\Gamma}_p^{-1} \mathbf{N}/[2(\sigma^2 \tilde{L} - \tilde{K})]$ minimizes the limit $S(\cdot)$ in Lemma A.5, the result follows by Remark 1 in [10].  □

REMARK 1.  $R$-estimators of linear regression parameters are also consistent and asymptotically normal [14]. Note, however, that the conditions placed on $\lambda$ in assumption A4 are slightly stronger than those placed on the weight function in [14], where the weight function is square integrable, not necessarily bounded or continuous. The conditions in A4 can be relaxed to some extent at the expense of stronger assumptions on $f$, but we do not pursue those extensions here. Since piecewise continuous and unbounded weight functions on $(0, 1)$ can be well approximated by differentiable, bounded weight functions, from a practical perspective, assumption A4 is not overly restrictive.

REMARK 2.  Using the Cauchy–Schwarz inequality,

$$\begin{aligned}
\sigma^2 \tilde{J} - \tilde{K}^2 &= \sigma^2 \mathrm{E}\{\lambda^2(F(Z_1))\} - (\mathrm{E}\{Z_1\lambda(F(Z_1))\})^2 \\
&\geq \sigma^2 \mathrm{E}\{\lambda^2(F(Z_1))\} - \mathrm{E}\{Z_1^2\}\mathrm{E}\{\lambda^2(F(Z_1))\} = 0,
\end{aligned} \tag{3.2}$$

with equality in (3.2) if and only if $\lambda$ is proportional to $F^{-1}$, which is not possible since $F^{-1}(0) = -\infty$, $F^{-1}(1) = \infty$ and $\lambda$ is bounded on $(0, 1)$. Hence, $\sigma^2 \tilde{J} - \tilde{K}^2 > 0$. $\tilde{K} = \int_0^1 F^{-1}(s)\lambda(s)\,ds$ is also greater than zero because $F^{-1}$ and $\lambda$ are strictly increasing functions on $(0, 1)$ and $\lambda$ is odd about $1/2$. Without assumption A5, $\sigma^2 \tilde{L} - \tilde{K}$ is not necessarily greater than zero, however. If the density function $f$ is differentiable, using integration by parts, it can be shown that

$$\tilde{L} = \mathrm{E}\{f(Z_1)\lambda'(F(Z_1))\} = -\int_{-\infty}^{\infty} f'(s)\lambda(F(s))\,ds = -\int_0^1 \frac{f'(F^{-1}(s))}{f(F^{-1}(s))}\lambda(s)\,ds.$$

Therefore, if $Z_1 \sim N(0, \sigma^2)$, then

$$\sigma^2 \tilde{L} = -\sigma^2 \int_0^1 \frac{f'(F^{-1}(s))}{f(F^{-1}(s))} \lambda(s) \, ds = \int_0^1 F^{-1}(s)\lambda(s) \, ds = \tilde{K}$$

and so A5 does not hold if $Z_1$ is Gaussian.

Remark 3. The asymptotic covariance matrix for $\hat{\boldsymbol{\phi}}_R$ is a scalar multiple of $n^{-1}\sigma^2 \boldsymbol{\Gamma}_p^{-1}$, the asymptotic covariance matrix for Gaussian likelihood estimators of the parameters of the corresponding $p$th-order autoregressive process. The same property holds for LAD and ML estimators of all-pass model parameters, as shown in [6] and [2], respectively. The LAD estimators are quasi-maximum likelihood estimators which can be obtained by maximizing the log-likelihood of an all-pass model with Laplace noise $(f(s) = \exp(-\sqrt{2}|s|/\sigma)/(\sqrt{2}\sigma))$. The appropriate scalar multiple is

$$(3.3) \qquad \frac{\text{Var}\,|Z_1|}{2(2\sigma^2 f(0) - E|Z_1|)^2}$$

in the LAD case ([6] contains an error in the calculation of the asymptotic variance; see [3] for the correction) and

$$(3.4) \qquad \frac{1}{2}\left(\sigma^2 \int_{-\infty}^{\infty} \frac{(f'(s))^2}{f(s)} \, ds - 1\right)^{-1}$$

in the ML case, while the multiple in (3.1) for $R$-estimation is

$$(3.5) \qquad \frac{\sigma^2 \tilde{J} - \tilde{K}^2}{2(\sigma^2 \tilde{L} - \tilde{K})^2}.$$

Consequently, the asymptotic relative efficiency (ARE) for $R$ to LAD is obtained by dividing (3.3) by (3.5) and the ARE for $R$ to ML is obtained by dividing (3.4) by (3.5).

Remark 4. Consider the sequence of weight functions $\{\lambda_m\}$ such that $\lambda_m(s) = 2\pi^{-1} \arctan(m(s - 1/2))$. It is straightforward to show that $\lambda_m$ satisfies assumptions A1 and A4 for all $m > 0$. If $I\{\cdot\}$ denotes the indicator function and $\tilde{\mu} := \text{median}\{Z_1\}$, then A5 is satisfied for large $m$ when $2\sigma^2 f(\tilde{\mu}) > E\{-Z_1 I\{Z_1 < \tilde{\mu}\} + Z_1 I\{Z_1 > \tilde{\mu}\}\}$; this holds for many distributions, including the Laplace, logistic and Student's $t$- (with degrees of freedom greater than two) distributions, and various asymmetric distributions $[0.4N(-1, 1) + 0.6N(2/3, 3^2)$ is one example]. Because $\lambda_m(s)$ converges pointwise to $-I\{s < 1/2\} + I\{s > 1/2\}$ on $(0, 1)$ as $m \to \infty$, $\tilde{J}_m = \int_0^1 \lambda_m^2(s) \, ds \to 1$ and, if $Z_1$ has median zero, $\tilde{K}_m = E\{Z_1 \lambda_m(F(Z_1))\} \to E|Z_1|$ and $\tilde{L}_m = E\{f(Z_1)\lambda_m'(F(Z_1))\} \to 2f(0)$. Hence, if $Z_1$ has median zero,

$$\frac{\sigma^2 \tilde{J}_m - \tilde{K}_m^2}{2(\sigma^2 \tilde{L}_m - \tilde{K}_m)^2} \to \frac{\sigma^2 - E^2|Z_1|}{2(2\sigma^2 f(0) - E|Z_1|)^2} = \frac{\text{Var}\,|Z_1|}{2(2\sigma^2 f(0) - E|Z_1|)^2}$$

and so $R$-estimation has virtually the same asymptotic efficiency as LAD estimation when the weight function $\lambda_m$ is used with $m$ large. If $Z_1$ has a Laplace distribution, LAD estimation corresponds to ML estimation. In the case of Laplace noise, therefore, $R$-estimation with weight function $\lambda_m$ and $m$ large also has essentially the same asymptotic efficiency as ML estimation.

REMARK 5.   Under general conditions, it can be shown that (3.5) equals (3.4) when the weight function is proportional to $-f'(F^{-1}(s))/f(F^{-1}(s))$. Thus, $R$-estimation has the same asymptotic efficiency as ML estimation when an optimal weight function $\lambda_f(s) \propto -f'(F^{-1}(s))/f(F^{-1}(s))$ is used. $\lambda_f$ is also an optimal weight function in the case of $R$-estimation for linear regression parameters (see, e.g., [15]). Note that if $Z_1$ has a Laplace distribution, then $\lambda_f(s) \propto -I\{s < 1/2\} + I\{s > 1/2\}$ for $s \in (0, 1/2) \cup (1/2, 1)$ ($\lambda_f$ does not exist at $s = 1/2$).

If $Z_1$ has a logistic distribution, then $f(s) = \pi/(\sqrt{3}\sigma) \exp(-s\pi/(\sqrt{3}\sigma))/[1 + \exp(-s\pi/(\sqrt{3}\sigma))]^2$ and so an optimal weight function $\lambda_f$ is given by the Wilcoxon weight function $\lambda(s) = s - 1/2$. For the Wilcoxon weights, assumption A5 is satisfied when $\sigma^2 \mathrm{E}\{f(Z_1)\} > \mathrm{E}\{Z_1 F(Z_1)\}$, which holds for the Laplace, logistic, Student's $t$- and $0.4\mathrm{N}(-1, 1) + 0.6\mathrm{N}(2/3, 3^2)$ distributions, as well as many others. Columns 2 and 3 of Table 1 give values of ARE for $R$ (with Wilcoxon weights) to LAD and $R$ (with Wilcoxon weights) to ML for a number of distributions. For the logistic and Student's $t$-distributions, $R$-estimation is asymptotically much more efficient than LAD and essentially as efficient as ML. Also, even though ML estimation is asymptotically 40% more efficient than $R$-estimation (with Wilcoxon weights) when the noise distribution is Laplace, $R$-estimation can still be useful in this case because $D(\cdot)$ tends to be smoother than $\sum_{t=1}^{n-p} |z_t(\cdot)|$ and hence easier to minimize. Figure 1 shows ML and $R$ objective functions for a realization of length $n = 50$ from an all-pass model with $p = 1$, $\phi_{01} = 0.5$ and Laplace noise with variance one. Observe that the ML objective function has many local minima and thus could be difficult to minimize using numerical optimization techniques.

REMARK 6.   Another weight function commonly used for $R$-estimation is the van der Waerden weight function $\lambda(s) = \Phi^{-1}(s)$, where $\Phi$ is the standard normal distribution function. Using results in [12], it can be shown that, if $f$ is absolutely continuous and almost everywhere differentiable with $0 < \int_{-\infty}^{\infty} (f'(s))^2/f(s)\,ds < \infty$, then A5 holds for the van der Waerden weights if and only if $Z_1$ is non-Gaussian. So, although A4 does not hold because $\Phi^{-1}$ is unbounded on $(0, 1)$, a bounded weight function approximating $\Phi^{-1}$ which does satisfy the assumptions can be found for a large class of non-Gaussian noise distributions. However, since $\Phi^{-1}$ is optimal when $Z_1 \sim \mathrm{N}(0, \sigma^2)$ and the parameters of a Gaussian all-pass series are not identifiable, the van der Waerden weights are not particularly useful for all-pass parameter estimation. Column 4 of Table 1 gives the ARE's for $R$ (with Wilcoxon weights) to $R$ (with van der Waerden weights) for various noise distributions. The van der Waerden weights are asymptotically superior to the Wilcoxon weights only when the distribution is close to Gaussian.

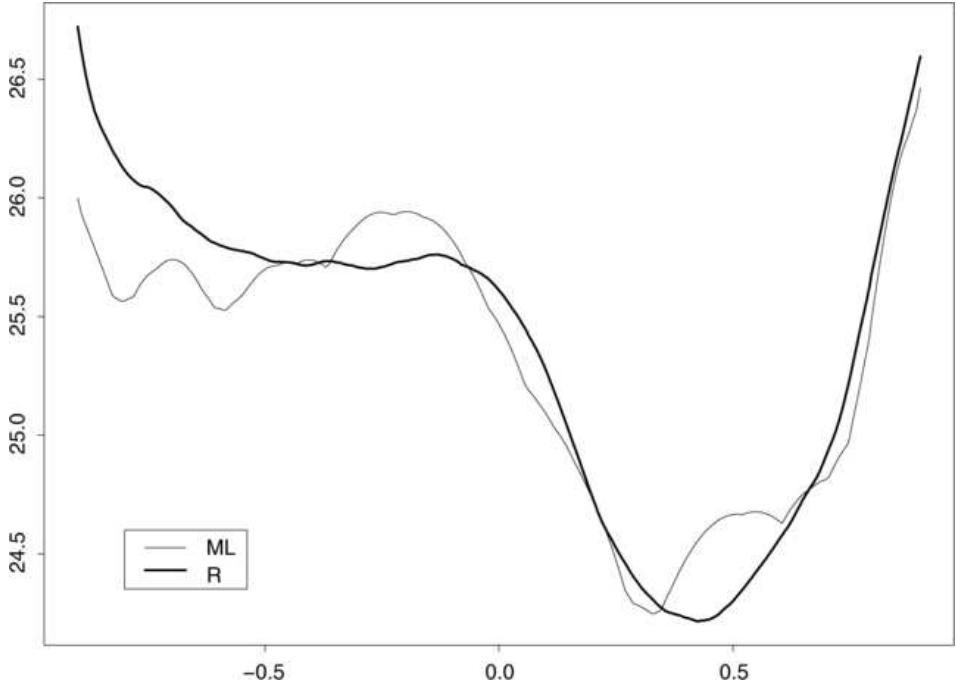

Fig. 1. *ML and R (with Wilcoxon weights) objective functions for a realization of length $n = 50$ from an all-pass model with $p = 1$, $\phi_{01} = 0.5$ and Laplace noise with variance one.*

3.2. *Order selection.* In practice, the true order $r$ of an all-pass model is usually unknown and must be estimated. In this section, we give an order selection procedure that is analogous to using the partial autocorrelation function to identify the order of an autoregressive model. First, note that

$$\frac{\sigma^2 \tilde{J} - \tilde{K}^2}{2(\sigma^2 \tilde{L} - \tilde{K})^2} = \frac{\tilde{J} - (|\phi_{0r}|/\sigma \tilde{K}_z)^2}{2(\sigma/|\phi_{0r}|\tilde{L}_z - |\phi_{0r}|/\sigma \tilde{K}_z)^2},$$

where $\tilde{K}_z := \int_0^1 F_z^{-1}(s)\lambda(s)\,ds$, $\tilde{L}_z := \int_0^1 f_z(F_z^{-1}(s))\lambda'(s)\,ds$ and $f_z$ and $F_z$ are the density and distribution functions, respectively, for $z_1 = \phi_{0r}^{-1} Z_1$. Be-cause $\hat{\phi}_R \xrightarrow{P} \phi_0$,

$$(3.6) \qquad \hat{s} := \left(\frac{1}{n}\sum_{t=1}^{n-p} z_t^2(\hat{\phi}_R)\right)^{1/2} \xrightarrow{P} (\mathrm{E}\{z_1^2\})^{1/2} = \frac{\sigma}{|\phi_{0r}|}$$

and $\hat{K}_z := n^{-1}D(\hat{\phi}_R) \xrightarrow{P} \tilde{K}_z$ by Lemma A.6 in the Appendix. Corollary 3.1 provides a consistent estimator of $\tilde{L}_z$.

COROLLARY 3.1. *Consider the kernel density estimator of $f_z$*

$$(3.7) \qquad \hat{f}_n(s) := \frac{1}{b_n n}\sum_{t=1}^{n-p} \kappa\left(\frac{s - z_t(\hat{\phi}_R)}{b_n}\right),$$

*where $\kappa$ is a uniformly continuous, differentiable kernel density function on $\mathbb{R}$ such that $\int |s\ln|s||^{1/2}|\kappa'(s)|\,ds < \infty$ and $\kappa'$ is uniformly continuous on $\mathbb{R}$, and where the bandwidth sequence $\{b_n\}$ is chosen so that $b_n \xrightarrow{P} 0$ and $b_n^2\sqrt{n} \xrightarrow{P} \infty$ as $n \to \infty$. If A1–A5 hold, then $\hat{L}_z := n^{-1}\sum_{t=1}^{n-p} \lambda'(t/(n-p))\hat{f}_n(z_{(t)}(\hat{\phi}_R)) \xrightarrow{P} \tilde{L}_z$.*

PROOF. If $\hat{F}_n(s) := (n-p)^{-1}\sum_{t=1}^{n-p} I\{z_t(\hat{\phi}_R) \le s\}$, $\hat{F}_n^{-1}(s) := \inf\{x: \hat{F}_n(x) \ge s\}$ and

$$\lambda'_n(s) := \lambda'\left(\frac{t}{n-p}\right) \qquad \text{for } s \in \left(\frac{t-1}{n-p}, \frac{t}{n-p}\right], \ t = 1, \ldots, n-p,$$

TABLE 1
*AREs for R (with Wilcoxon weights) to LAD, R (with Wilcoxon weights) to ML and R (with Wilcoxon weights) to R (with van der Waerden weights) for the Laplace distribution, the logistic distribution and Student's t-distribution with several different degrees of freedom*

| Noise distribution | ARE (R to LAD) | ARE (R to ML) | ARE (R to R) |
|---|---|---|---|
| Laplace | 0.600 | 0.600 | 1.026 |
| logistic | 1.976 | 1.000 | 1.049 |
| $t(3)$ | 1.411 | 0.962 | 1.208 |
| $t(6)$ | 2.068 | 0.997 | 1.083 |
| $t(9)$ | 2.354 | 0.980 | 1.023 |
| $t(12)$ | 2.510 | 0.964 | 0.990 |
| $t(15)$ | 2.607 | 0.952 | 0.971 |
| $t(20)$ | 2.707 | 0.937 | 0.953 |
| $t(30)$ | 2.810 | 0.921 | 0.938 |

then $n\hat{L}_z/(n-p) = \int_0^1 \hat{f}_n(\hat{F}_n^{-1}(s))\lambda_n'(s)\,ds$. By the uniform continuity of $\lambda'$, $\sup_{s\in(0,1)} |\lambda_n'(s) - \lambda'(s)| \to 0$. Consequently, since $\sup_{s\in(0,1)} f_z(F_z^{-1}(s)) < \infty$ and $\sup_{s\in(0,1)} |\lambda'(s)| < \infty$, the proof is complete if

$$\sup_{s\in(0,1)} |\hat{f}_n(\hat{F}_n^{-1}(s)) - f_z(F_z^{-1}(s))|$$

$$(3.8) \qquad \leq \sup_{s\in(0,1)} |\hat{f}_n(\hat{F}_n^{-1}(s)) - f_z(\hat{F}_n^{-1}(s))|$$

$$+ \sup_{s\in(0,1)} |f_z(\hat{F}_n^{-1}(s)) - f_z(F_z^{-1}(s))|$$

is $o_p(1)$. Because $\sup_{s\in\mathbb{R}} |\hat{f}_n(s) - f_z(s)| \xrightarrow{P} 0$ (for proof of this result, see Lemma 16 on page 88 of [3]; a similar result is given in Theorem 3 of [21]), the first term in (3.8) is $o_p(1)$. We now consider the second term and use an argument similar to one found in the proof of Lemma 4 in [18]. Note that $\sup_{s\in(0,1)} |F_z(\hat{F}_n^{-1}(s)) - s| = \sup_{s\in\mathbb{R}} |\hat{F}_n(s) - F_z(s)|$ and, using the Glivenko–Cantelli theorem, it can be shown that $\sup_{s\in\mathbb{R}} |\hat{F}_n(s) - F_z(s)| \xrightarrow{P} 0$. Therefore, because $f_z(F_z^{-1}(\cdot))$ is uniformly continuous on $(0,1)$ and $F_z^{-1}(F_z(s)) = s$ for all $s \in \mathbb{R}$ (since $F_z$ is strictly increasing on $\mathbb{R}$), we have

$$\sup_{s\in(0,1)} |f_z(\hat{F}_n^{-1}(s)) - f_z(F_z^{-1}(s))|$$

$$= \sup_{s\in(0,1)} |f_z(F_z^{-1}[F_z\{\hat{F}_n^{-1}(s)\}]) - f_z(F_z^{-1}(s))| \xrightarrow{P} 0. \qquad \square$$

It follows that

$$(3.9) \qquad \frac{\tilde{J} - (\hat{s}^{-1}\hat{K}_z)^2}{2(\hat{s}\hat{L}_z - \hat{s}^{-1}\hat{K}_z)^2} \xrightarrow{P} \frac{\sigma^2\tilde{J} - \tilde{K}^2}{2(\sigma^2\tilde{L} - \tilde{K})^2}.$$

Note that the Gaussian and the Student $t$-densities satisfy the conditions for the kernel density function $\kappa$ in Corollary 3.1.

We now give the following corollary for use in order selection.

COROLLARY 3.2. *Assume* A1–A5 *hold. If the true order of the all-pass model is $r$ and the order of the fitted model is $p > r$, then $n^{1/2}\hat{\phi}_{p,R} \xrightarrow{d}$ $\mathrm{N}(0, (\sigma^2\tilde{J} - \tilde{K}^2)/[2(\sigma^2\tilde{L} - \tilde{K})^2])$.*

PROOF. By Problem 8.15 in [7], the $p$th diagonal element of $\mathbf{\Gamma}_p^{-1}$ is $\sigma^{-2}$ if $p > r$, so the result follows from (3.1). $\square$

A practical approach to order determination using a large sample is described as follows:

TABLE 2
*Empirical means, standard deviations and percent coverages of
nominal 95% confidence intervals for R-estimates of all-pass model parameters.
The LAD-like score function $\lambda(s) = 2\pi^{-1}\arctan(500(s-1/2))$ and the Wilcoxon
score function $\lambda(s) = s - 1/2$ were used. The noise distribution is Laplace with variance
one*

| | | Asymptotic | | Empirical | |
|---|---|---|---|---|---|
| $n$ | mean | std. dev. (LAD/Wilcoxon) | mean (LAD/Wilcoxon) | std. dev. (LAD/Wilcoxon) | % coverage (LAD/Wilcoxon) |
| 500 | $\phi_1 = 0.5$ | 0.0275/0.0354 | 0.499/0.497 | 0.0332/0.0593 | 97.7/96.2 |
| 5000 | $\phi_1 = 0.5$ | 0.0087/0.0112 | 0.500/0.499 | 0.0093/0.0112 | 97.9/96.0 |
| 500 | $\phi_1 = 0.3$ | 0.0291/0.0374 | 0.299/0.299 | 0.0413/0.0444 | 96.5/94.9 |
| | $\phi_2 = 0.4$ | 0.0291/0.0374 | 0.397/0.392 | 0.0479/0.0599 | 97.6/95.4 |
| 5000 | $\phi_1 = 0.3$ | 0.0092/0.0118 | 0.300/0.300 | 0.0101/0.0122 | 97.6/95.2 |
| | $\phi_2 = 0.4$ | 0.0092/0.0118 | 0.399/0.399 | 0.0099/0.0119 | 97.5/96.7 |

1. For some large $P$, fit all-pass models of order $p$, $p = 1, 2, \ldots, P$, via R-estimation and obtain the $p$th coefficient, $\hat{\phi}_{p,R}$, for each.
2. Let the model order $r$ be the smallest order beyond which the estimated coefficients are statistically insignificant; that is, $r = \min\{0 \le p \le P : |\hat{\phi}_{j,R}| < 1.96\hat{\tau}n^{-1/2}$ for $j > p\}$, where $\hat{\tau} := ([\tilde{J} - (\hat{s}^{-1}\hat{K}_z)^2]/[2(\hat{s}\hat{L}_z - \hat{s}^{-1}\hat{K}_z)^2])^{1/2}$ and the estimates $\hat{s}$, $\hat{K}_z$ and $\hat{L}_z$ are from the fitted $P$th-order model.

## 4. Numerical results.

4.1. *Simulation study.* In this section we give the results of a simulation study to assess the quality of the asymptotic approximations for finite samples. First, for each of 1000 replicates, we simulated all-pass data and found $\hat{\phi}_R$ by minimizing $D$ in (2.4). To reduce the possibility of the optimizer getting trapped at local minima, we chose 1000 random starting values for each replicate. We evaluated $D$ at each of the 1000 candidate values and then reduced the collection of initial values to the twelve with the smallest values of $D$. Optimized values were found using these twelve initial values as starting points. The optimized value for which $D$ was smallest was chosen to be $\hat{\phi}_R$. Confidence intervals for the elements of $\phi_0$ were constructed using (3.1) and the estimator in (3.9). For the kernel density estimator (3.7), we used the standard Gaussian kernel density function and, because of its recommendation in [22], page 48, we used bandwidth $b_n = 0.9n^{-1/5}\min\{\hat{s}, IQR/1.34\}$, where $\hat{s}$, defined in (3.6), is the sample standard deviation for $\{z_t(\hat{\phi}_R)\}$ and $IQR$ is the interquartile range for $\{z_t(\hat{\phi}_R)\}$.

Results of these simulations appear in Tables 2 and 3. We show the empirical means, standard deviations, and percent coverages of nominal 95% confidence intervals for the $R$-estimates of all-pass model parameters. The LAD-like score function $\lambda(s) = 2\pi^{-1}\arctan(500(s-1/2))$ and the Wilcoxon score function $\lambda(s) = s - 1/2$ were used. Asymptotic means and standard deviations were obtained using Theorem 3.1. Note that the $R$-estimates appear nearly unbiased and the confidence interval coverages are close to the nominal 95% level. The asymptotic standard deviations tend to understate the true variability of the estimates when $n = 500$, but are fairly accurate when $n = 5000$. Normal probability plots show that the $R$-estimates are approximately normal, particularly when $n = 5000$. The quality of the asymptotic approximations for finite samples is similar for LAD and ML estimates (see [2] and [6]).

We also ran simulations to assess the order selection procedure described in Section 3.2. For each of 100 replicates, we simulated all-pass data and

TABLE 3

*Empirical means, standard deviations and percent coverages of nominal 95% confidence intervals for R-estimates of all-pass model parameters. The LAD-like score function $\lambda(s) = 2\pi^{-1}\arctan(500(s-1/2))$ and the Wilcoxon score function $\lambda(s) = s - 1/2$ were used. The noise distribution is Student's t with three degrees of freedom*

| | Asymptotic | | Empirical | | |
|---|---|---|---|---|---|
| $n$ | mean | std. dev. (LAD/Wilcoxon) | mean (LAD/Wilcoxon) | std. dev. (LAD/Wilcoxon) | % coverage (LAD/Wilcoxon) |
| 500 | $\phi_1 = 0.5$ | 0.0327/0.0279 | 0.499/0.498 | 0.0405/0.0331 | 95.8/96.2 |
| 5000 | $\phi_1 = 0.5$ | 0.0103/0.0088 | 0.500/0.500 | 0.0110/0.0090 | 95.2/95.6 |
| 500 | $\phi_1 = 0.3$ | 0.0346/0.0296 | 0.301/0.299 | 0.0403/0.0366 | 95.1/94.7 |
| | $\phi_2 = 0.4$ | 0.0346/0.0296 | 0.396/0.396 | 0.0418/0.0366 | 95.2/94.9 |
| 5000 | $\phi_1 = 0.3$ | 0.0109/0.0093 | 0.300/0.300 | 0.0118/0.0095 | 94.0/95.4 |
| | $\phi_2 = 0.4$ | 0.0109/0.0093 | 0.400/0.400 | 0.0115/0.0097 | 94.6/95.1 |

TABLE 4

*The frequencies for each estimate of model order $r$ when $P = 5$ and the Wilcoxon scores were used*

| | | Laplace noise | | | | | | $t$ noise | | | | | |
|---|---|---|---|---|---|---|---|---|---|---|---|---|---|
| $n$ | Model parameters | 0 | 1 | 2 | 3 | 4 | 5 | 0 | 1 | 2 | 3 | 4 | 5 |
| 500 | $\phi_0 = 0.5$ | 0 | 58 | 7 | 10 | 8 | 17 | 0 | 52 | 8 | 9 | 13 | 18 |
| 5000 | $(r = 1)$ | 0 | 67 | 2 | 3 | 7 | 21 | 0 | 56 | 0 | 3 | 7 | 34 |
| 500 | $\phi_0 = (0.3, 0.4)'$ | 0 | 0 | 69 | 1 | 16 | 14 | 0 | 0 | 57 | 13 | 16 | 14 |
| 5000 | $(r = 2)$ | 0 | 0 | 82 | 5 | 5 | 8 | 0 | 0 | 81 | 5 | 6 | 8 |

estimated the model order $r$ using the procedure in Section 3.2 with $P = 5$ and Wilcoxon scores. Table 4 gives the frequencies for each estimate of $r$. In all cases, the procedure appears to be fairly successful at identifying the true value of $r$. Model orders less than $r$ were never selected, so underestimating $r$ is clearly not a concern.

4.2. *Deconvolution.* Applications for all-pass models are not limited to uncorrelated time series. As discussed in Section 1 and in [2], all-pass models can also be used to identify and model noncausal and noninvertible ARMA series. If, for example, a causal, invertible ARMA model is fit to a causal, noninvertible series, the residuals follow a causal all-pass model of order $r$, where $r$ is the number of roots of the true moving average polynomial inside the unit circle. Therefore, the order of noninvertibility of the ARMA, $r$, can be determined by identifying the all-pass order of the residuals.

Consider the simulated water gun seismogram $\{X_t\}_{t=1}^{1000}$ shown in Figure 2(a), where $X_t = \sum_k \beta_k Z_{t-k}$, $\{\beta_k\}$ is the water gun wavelet sequence in Figure 8(2) of [19] and $\{Z_t\}$ is a reflectivity sequence which was simulated as i.i.d. noise from the Student $t$-distribution with five degrees of freedom. Andrews, Davis and Breidt [2] modeled $\{X_t\}$ as a possibly noninvertible ARMA, using ML estimation for all-pass models to identify an appropriate order of noninvertibility. The wavelet and reflectivity sequences were then reconstructed from $\{X_t\}$ using the fitted ARMA model. This deconvolution procedure is of interest because, for an observed water gun seismogram, the reflectivity sequence is unknown and corresponds to reflection coefficients for layers of the Earth. In this section, we identify an appropriate order of noninvertibility for $\{X_t\}$ using $R$-estimation for all-pass models and we compare the $R$-estimation results to ML results in [2].

Andrews, Davis and Breidt [2] first fit a causal, invertible ARMA$(12, 13)$ model $\phi(B)X_t = \theta(B)W_t$ to the simulated seismogram $\{X_t\}$ using Gaussian ML. The residuals from this fitted ARMA model are denoted $\{\hat{W}_t\}$. From the sample autocorrelation functions for $\{\hat{W}_t\}$, $\{\hat{W}_t^2\}$ and $\{|\hat{W}_t|\}$ in Figure 2(b)–(d), it appears that these ARMA residuals are uncorrelated but dependent, suggesting that a causal, invertible model is inappropriate for $\{X_t\}$. Using ML estimation and the Student $t$-density, a causal all-pass model of order two was determined to be most suitable for $\{\hat{W}_t\}$ [2]. The ML estimates of the all-pass model parameters are $\hat{\phi}_{ML} = (1.5286, -0.5908)'$, both with standard error 0.0338. Since the all-pass residuals appear independent, Andrews, Davis and Breidt [2] concluded that a causal, noninvertible ARMA(12,13) with two roots of the moving average polynomial inside the unit circle is an appropriate model for $\{X_t\}$.

When the Wilcoxon weight function, the standard Gaussian kernel density function and bandwidth $b_n = 0.9n^{-1/5}\min\{\hat{s}, IQR/1.34\}$ are used, the order selection procedure described in Section 3.2 also indicates that an all-pass model of order two is appropriate for $\{\hat{W}_t\}$. The $R$-estimates of the all-pass model parameters are $\hat{\boldsymbol{\phi}}_R = (1.5052, -0.5700)'$, both with standard error 0.0343. In Figure 3, we show the sample autocorrelation functions for the squares and absolute values of $\{\hat{Z}_t\}$, the residuals from the all-pass model fit to $\{\hat{W}_t\}$ using $R$-estimation; these all-pass residuals appear independent. Therefore, in this example, the ML and $R$ all-pass estimation results are nearly identical, even though no specific distributional information was used for $R$-estimation.

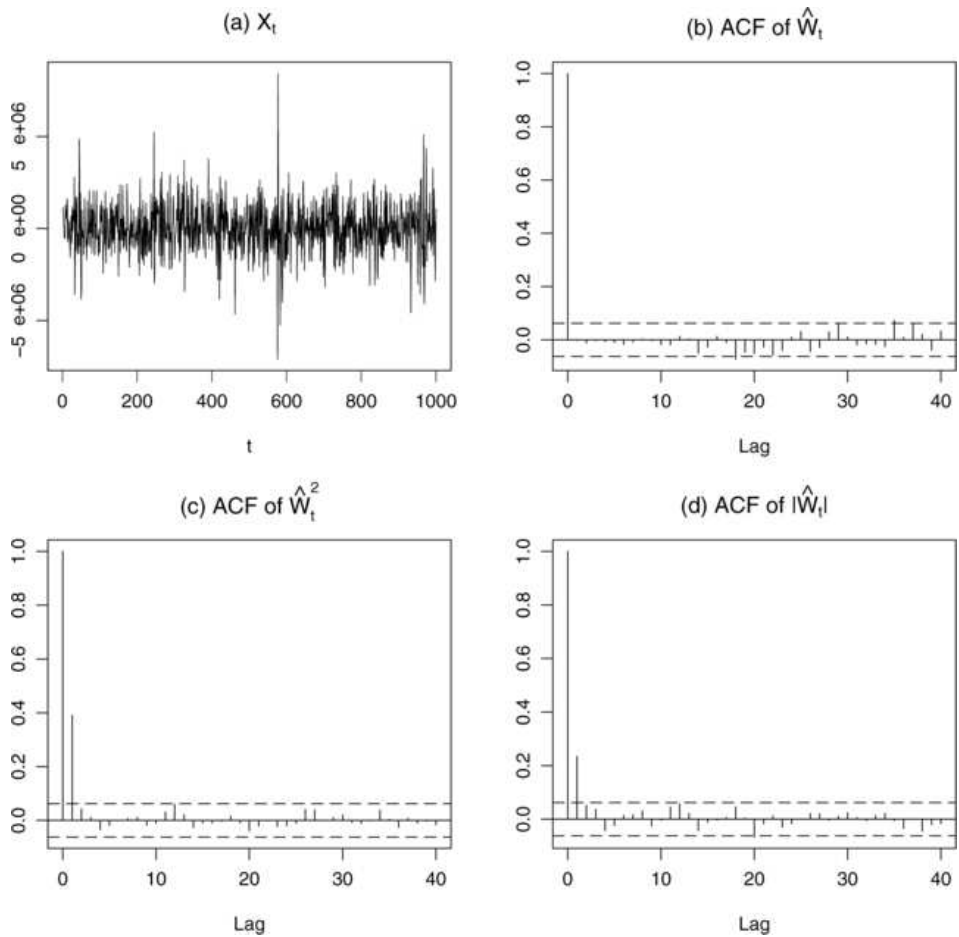

Fig. 2. (a) *The simulated seismogram of length* 1000, $\{X_t\}$, *and the sample autocorrelation functions with bounds* $\pm 1.96/\sqrt{1000}$ *for* (b) $\{\hat{W}_t\}$, (c) $\{\hat{W}_t^2\}$ *and* (d) $\{|\hat{W}_t|\}$.

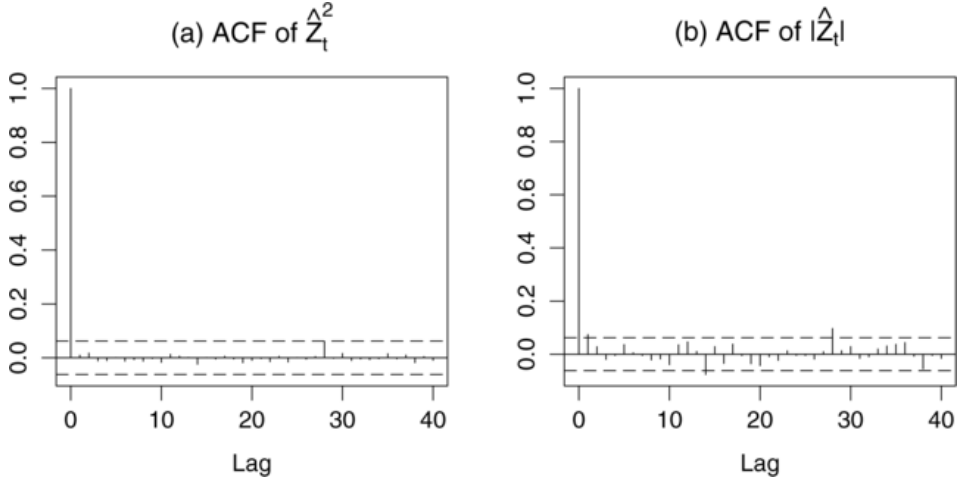

Fig. 3. *Diagnostics for the all-pass model of order two fit to the causal, invertible ARMA residuals using R-estimation. The sample autocorrelation functions with bounds $\pm 1.96/\sqrt{1000}$ for* (a) $\{\hat{Z}_t^2\}$ *and* (b) $\{|\tilde{Z}_t|\}$ *are shown.*

## APPENDIX

This section contains proofs of the lemmas used to establish the results of Section 3. We assume that assumptions A1–A5 hold throughout. First, note that for $j \in \{1, \ldots, p\}$ and $t \in \{1, \ldots, n-p\}$,

$$\text{(A.1)} \qquad \frac{\partial z_t(\phi)}{\partial \phi_j} = \frac{1}{\phi(B^{-1})}\{X_{t+p-j} + z_{t+j}(\phi)\}$$

(see [2] for details). Evaluating (A.1) at the true value of $\phi$ and ignoring the effect of recursion initialization, we have

$$\frac{\partial z_t(\phi_0)}{\partial \phi_j} = \frac{1}{\phi_0(B^{-1})}\left\{\frac{-\phi_0(B^{-1})B^p z_{t+p-j}}{\phi_0(B)} + z_{t+j}(\phi_0)\right\}$$

$$\text{(A.2)} \qquad \simeq \frac{-z_{t-j}}{\phi_0(B)} + \frac{z_{t+j}}{\phi_0(B^{-1})},$$

where the first term is an element of $\sigma(z_{t-1}, z_{t-2}, \ldots)$ and the second term is an element of $\sigma(z_{t+1}, z_{t+2}, \ldots)$ because $\phi_0(B)$ is a causal operator and $\phi_0(B^{-1})$ is a purely noncausal operator. It follows that (A.2) is independent of $z_t = \phi_{0r}^{-1} Z_t$. Thus, if $F_z$ is the distribution function of $z_1$ and $g_t(\phi) := \lambda(F_z(z_t))z_t(\phi)$, then for $j \in \{1, \ldots, p\}$,

$$\frac{\partial g_t(\phi_0)}{\partial \phi_j} = \lambda(F_z(z_t))\frac{\partial z_t(\phi_0)}{\partial \phi_j} \simeq \lambda(F_z(z_t))\left\{\frac{-z_{t-j}}{\phi_0(B)} + \frac{z_{t+j}}{\phi_0(B^{-1})}\right\} =: \frac{\partial g_t^*(\phi_0)}{\partial \phi_j}.$$

The expected value of $\partial g_t^*(\phi_0)/\partial \phi_j$ is zero by the independence of its two terms.

We now compute the autocovariance function $\gamma^\dagger(h)$ of the zero-mean, stationary process $\{\mathbf{u}'\partial g_t^*(\boldsymbol{\phi}_0)/\partial\boldsymbol{\phi}\}$ for $\mathbf{u} \in \mathbb{R}^p$:

$$\gamma^\dagger(h) = \mathrm{E}\left\{\mathbf{u}'\frac{\partial g_t^*(\boldsymbol{\phi}_0)}{\partial\boldsymbol{\phi}}\left(\frac{\partial g_{t+h}^*(\boldsymbol{\phi}_0)}{\partial\boldsymbol{\phi}}\right)'\mathbf{u}\right\} = \mathbf{u}'[\nu_{jk}(h)]_{j,k=1}^p\mathbf{u},$$

where

$$\nu_{jk}(h) := \begin{cases} 2\phi_{0r}^{-2}\tilde{J}\gamma(j-k), & \text{if } h = 0, \\ -\psi_{|h|-j}\psi_{|h|-k}\phi_{0r}^{-2}\tilde{K}^2, & \text{if } h \neq 0, \end{cases}$$

and the $\psi_l$ are given by $\sum_{l=0}^\infty \psi_l z^l = 1/\phi_0(z)$ with $\psi_l = 0$ for $l < 0$. Thus,

$$\gamma^\dagger(0) + 2\sum_{h=1}^\infty \gamma^\dagger(h)$$

$$= \mathbf{u}'\left\{[2\phi_{0r}^{-2}\tilde{J}\gamma(j-k)]_{j,k=1}^p - 2\phi_{0r}^{-2}\tilde{K}^2\left[\sum_{h=1}^\infty \psi_{h-j}\psi_{h-k}\right]_{j,k=1}^p\right\}\mathbf{u}$$

$$= 2\phi_{0r}^{-2}(\sigma^2\tilde{J} - \tilde{K}^2)\mathbf{u}'\sigma^{-2}\boldsymbol{\Gamma}_p\mathbf{u}.$$

The preceding calculations lead directly to the following lemma.

LEMMA A.1. *As* $n \to \infty$, $n^{-1/2}\sum_{t=1}^{n-p}\partial g_t(\boldsymbol{\phi}_0)/\partial\boldsymbol{\phi} \xrightarrow{d} \mathbf{N} \sim \mathrm{N}(\mathbf{0}, 2\phi_{0r}^{-2} \times \{\sigma^2\tilde{J} - \tilde{K}^2\}\sigma^{-2}\boldsymbol{\Gamma}_p)$.

PROOF. Note that, for $t \in \{0, \ldots, n-p-1\}$,

$$z_{n-p-t} = \sum_{l=0}^\infty \psi_l(\phi_0(B^{-1})z_{n-p-t+l} \quad \text{and}$$

(A.3)

$$z_{n-p-t}(\boldsymbol{\phi}_0) = \sum_{l=0}^t \psi_l(\phi_0(B^{-1})z_{n-p-t+l}).$$

Because there exist constants $c > 0$ and $0 < d < 1$ such that $|\psi_l| < cd^l$ for all $l \in \{0, 1, \ldots\}$ (see [7], Section 3.3), we have

$$\sum_{t=1}^{n-p}\mathrm{E}\left|\frac{\partial g_t(\boldsymbol{\phi}_0)}{\partial\phi_j} - \frac{\partial g_t^*(\boldsymbol{\phi}_0)}{\partial\phi_j}\right| = \sum_{t=1}^{n-p}\mathrm{E}\left|\lambda(F_z(z_t))\left\{\frac{z_{t+j}(\boldsymbol{\phi}_0)}{\phi_0(B^{-1})} - \frac{z_{t+j}}{\phi_0(B^{-1})}\right\}\right| = O(1)$$

for $j \in \{1, \ldots, p\}$. Consequently, $n^{-1/2}\sum_{t=1}^{n-p}[\partial g_t(\boldsymbol{\phi}_0)/\partial\boldsymbol{\phi} - \partial g_t^*(\boldsymbol{\phi}_0)/\partial\boldsymbol{\phi}] \to \mathbf{0}$ in $L_1$ and hence in probability.

Let $\mathbf{u} \in \mathbb{R}^p$. By the Cramér–Wold device, it suffices to show that $n^{-1/2} \times \sum_{t=1}^{n-p}\mathbf{u}'\partial g_t^*(\boldsymbol{\phi}_0)/\partial\boldsymbol{\phi} \xrightarrow{d} \mathbf{u}'\mathbf{N} \sim \mathrm{N}(0, 2\phi_{0r}^{-2}\{\sigma^2\tilde{J} - \tilde{K}^2\}\mathbf{u}'\sigma^{-2}\boldsymbol{\Gamma}_p\mathbf{u})$. Elements of the infinite-order moving average stationary sequence $\{\mathbf{u}'\partial g_t^*(\boldsymbol{\phi}_0)/\partial\boldsymbol{\phi}$ can be truncated to create a finite-order moving average stationary sequence. By applying a central limit theorem ([7], Theorem 6.4.2) to each truncation level, asymptotic normality can be deduced. The details are omitted.  □

Now, consider the mixed partials of $g_t(\boldsymbol{\phi})$. For $j, k \in \{1, \ldots, p\}$,

$$\frac{\partial^2 z_t(\boldsymbol{\phi}_0)}{\partial \phi_j \, \partial \phi_k} = \frac{1}{\phi_0^2(B^{-1})} \{X_{t+p+j-k} + X_{t+p+k-j} + 2z_{t+j+k}(\boldsymbol{\phi}_0)\}$$

$$\simeq \frac{-z_{t+j-k} - z_{t+k-j}}{\phi_0(B^{-1})\phi_0(B)} + \frac{2z_{t+j+k}}{\phi_0^2(B^{-1})}$$

$$= -\sum_{m=0}^{\infty} \sum_{\ell=0}^{\infty} \psi_m \psi_\ell (z_{t+j-k-\ell+m} + z_{t+k-j-\ell+m}) + \frac{2z_{t+j+k}}{\phi_0^2(B^{-1})},$$

and so

$$\frac{\partial^2 g_t(\boldsymbol{\phi}_0)}{\partial \phi_j \, \partial \phi_k} = \lambda(F_z(z_t)) \frac{\partial^2 z_t(\boldsymbol{\phi}_0)}{\partial \phi_j \, \partial \phi_k}$$

$$\simeq \lambda(F_z(z_t)) \left\{ -\sum_{m=0}^{\infty} \sum_{\ell=0}^{\infty} \psi_m \psi_\ell (z_{t+j-k-\ell+m} + z_{t+k-j-\ell+m}) \right.$$

(A.4)
$$\left. + \frac{2z_{t+j+k}}{\phi_0^2(B^{-1})} \right\}$$

$$=: \frac{\partial^2 g_t^*(\boldsymbol{\phi}_0)}{\partial \phi_j \, \partial \phi_k}.$$

(A.4) has expectation $-2\sigma^{-2}\gamma(j-k) \int_0^1 F_z^{-1}(s)\lambda(s) \, ds = -2|\phi_{0r}|^{-1}\tilde{K}\sigma^{-2} \times \gamma(j-k)$.

LEMMA A.2. *As* $n \to \infty$, $n^{-1} \sum_{t=1}^{n-p} \partial^2 g_t(\boldsymbol{\phi}_0)/(\partial \boldsymbol{\phi} \, \partial \boldsymbol{\phi}') \xrightarrow{P} -2|\phi_{0r}|^{-1} \times \tilde{K}\sigma^{-2}\boldsymbol{\Gamma}_p$.

PROOF. It can be shown that $n^{-1} \sum_{t=1}^{n-p} [\partial^2 g_t(\boldsymbol{\phi}_0)/(\partial \boldsymbol{\phi} \, \partial \boldsymbol{\phi}') - \partial^2 g_t^*(\boldsymbol{\phi}_0)/(\partial \boldsymbol{\phi} \, \partial \boldsymbol{\phi}')] \to \mathbf{0}$ in $L_1$ and in probability. Because (A.4) has expectation $-2|\phi_{0r}|^{-1}\tilde{K}\sigma^{-2}\gamma(j-k)$, $n^{-1} \sum_{t=1}^{n-p} \partial^2 g_t^*(\boldsymbol{\phi}_0)/(\partial \boldsymbol{\phi} \partial \boldsymbol{\phi}') \xrightarrow{P} -2|\phi_{0r}|^{-1}\tilde{K}\sigma^{-2}\boldsymbol{\Gamma}_p$ by the ergodic theorem.  □

LEMMA A.3. *For any* $T \in (0, \infty)$, *as* $n \to \infty$,

(A.5)
$$\sup_{\|\mathbf{u}\| \le T} \left| n^{-1/2} \sum_{t=1}^{n-p} \mathbf{u}' \left[ \lambda\left( \frac{R_t(\boldsymbol{\phi}_0 + n^{-1/2}\mathbf{u})}{n-p+1} \right) - \lambda(F_z(z_t)) \right] \frac{\partial z_t(\boldsymbol{\phi}_0)}{\partial \boldsymbol{\phi}} \right.$$

$$\left. - 2|\phi_{0r}|^{-1}\tilde{L}\mathbf{u}'\boldsymbol{\Gamma}_p\mathbf{u} \right| \xrightarrow{P} 0.$$

PROOF.   Observe that the left-hand side of (A.5) is bounded above by

$$
\text{(A.6)} \quad \sup_{\|\mathbf{u}\| \le T} \left| \frac{1}{\sqrt{n}} \sum_{t=1}^{n-p} \mathbf{u}' \lambda'(F_z(z_t)) \frac{\partial z_t(\phi_0)}{\partial \phi} \left[ \frac{R_t(\phi_0 + n^{-1/2}\mathbf{u})}{n-p+1} - F_z(z_t) \right] \right.
$$

$$
\left. - 2|\phi_{0r}|^{-1} \tilde{L} \mathbf{u}' \mathbf{\Gamma}_p \mathbf{u} \right|
$$

$$
\text{(A.7)} \quad + \sup_{\|\mathbf{u}\| \le T} \left| \frac{1}{\sqrt{n}} \sum_{t=1}^{n-p} \mathbf{u}' [\lambda'(F_{t,n}^*(\mathbf{u})) - \lambda'(F_z(z_t))] \frac{\partial z_t(\phi_0)}{\partial \phi} \right.
$$

$$
\left. \times \left[ \frac{R_t(\phi_0 + n^{-1/2}\mathbf{u})}{n-p+1} - F_z(z_t) \right] \right|,
$$

where $F_{t,n}^*(\mathbf{u})$ is between $F_z(z_t)$ and $R_t(\phi_0 + n^{-1/2}\mathbf{u})/(n-p+1)$. If $F_n(x) := n^{-1} \sum_{t=1}^{n-p} I\{z_t \le x\}$, an upper bound for (A.6) is

$$
\text{(A.8)} \quad \sup_{\|\mathbf{u}\| \le T} \left| \frac{1}{\sqrt{n}} \sum_{t=1}^{n-p} \mathbf{u}' \lambda'(F_z(z_t)) \frac{\partial z_t(\phi_0)}{\partial \phi} \right.
$$

$$
\left. \times \left[ \frac{R_t(\phi_0 + n^{-1/2}\mathbf{u})}{n-p+1} - F_n\left( z_t\left( \phi_0 + \frac{\mathbf{u}}{\sqrt{n}} \right) \right) \right] \right|
$$

$$
\text{(A.9)} \quad + \sup_{\|\mathbf{u}\| \le T} \left| \frac{1}{\sqrt{n}} \sum_{t=1}^{n-p} \mathbf{u}' \lambda'(F_z(z_t)) \frac{\partial z_t(\phi_0)}{\partial \phi} \right.
$$

$$
\left. \times \left[ F_n\left( z_t\left( \phi_0 + \frac{\mathbf{u}}{\sqrt{n}} \right) \right) - F_z\left( z_t\left( \phi_0 + \frac{\mathbf{u}}{\sqrt{n}} \right) \right) \right] \right|
$$

$$
\text{(A.10)} \quad + \sup_{\|\mathbf{u}\| \le T} \left| \frac{1}{\sqrt{n}} \sum_{t=1}^{n-p} \mathbf{u}' \lambda'(F_z(z_t)) \frac{\partial z_t(\phi_0)}{\partial \phi} \right.
$$

$$
\times \left[ F_z\left( z_t\left( \phi_0 + \frac{\mathbf{u}}{\sqrt{n}} \right) \right) - F_z(z_t) \right]
$$

$$
\left. - 2|\phi_{0r}|^{-1} \tilde{L} \mathbf{u}' \mathbf{\Gamma}_p \mathbf{u} \right|.
$$

Because

$$
\sup_{\|\mathbf{u}\| \le T} \frac{1}{n} \sum_{t=1}^{n-p} \left| \mathbf{u}' \lambda'(F_z(z_t)) \frac{\partial z_t(\phi_0)}{\partial \phi} \right| = O_p(1)
$$

and, by Lemma 3 on page 55 of [3],

$$(A.11) \quad \sup_{\|\mathbf{u}\| \leq T, t \in \{1,\dots,n-p\}} \sqrt{n} \left| \frac{R_t(\boldsymbol{\phi}_0 + n^{-1/2}\mathbf{u})}{n-p+1} - F_n\left(z_t\left(\boldsymbol{\phi}_0 + \frac{\mathbf{u}}{\sqrt{n}}\right)\right) \right| \xrightarrow{P} 0,$$

(A.8) is $o_p(1)$. Lemma 10 on page 76 of [3] establishes that (A.9) is $o_p(1)$. Finally,

$$\sup_{\|\mathbf{u}\| \leq T} \left| \frac{1}{\sqrt{n}} \sum_{t=1}^{n-p} \mathbf{u}' \lambda'(F_z(z_t)) \frac{\partial z_t(\boldsymbol{\phi}_0)}{\partial \boldsymbol{\phi}} \left[ F_z\left(z_t\left(\boldsymbol{\phi}_0 + \frac{\mathbf{u}}{\sqrt{n}}\right)\right) - F_z(z_t) \right] \right.$$
$$\left. - 2|\phi_{0r}|^{-1}\tilde{L}\mathbf{u}'\boldsymbol{\Gamma}_p\mathbf{u} \right|$$

$$= \sup_{\|\mathbf{u}\| \leq T} \left| \frac{1}{\sqrt{n}} \sum_{t=1}^{n-p} \mathbf{u}' \lambda'(F_z(z_t)) \frac{\partial z_t(\boldsymbol{\phi}_0)}{\partial \boldsymbol{\phi}} f_z(z_{t,n}^*(\mathbf{u})) \left( z_t\left(\boldsymbol{\phi}_0 + \frac{\mathbf{u}}{\sqrt{n}}\right) - z_t \right) \right.$$
$$\left. - 2|\phi_{0r}|^{-1}\tilde{L}\mathbf{u}'\boldsymbol{\Gamma}_p\mathbf{u} \right|$$

$$\leq \sup_{\|\mathbf{u}\| \leq T} \left| \frac{1}{\sqrt{n}} \sum_{t=1}^{n-p} \mathbf{u}' \lambda'(F_z(z_t)) \frac{\partial z_t(\boldsymbol{\phi}_0)}{\partial \boldsymbol{\phi}} f_z(z_{t,n}^*(\mathbf{u}))(z_t(\boldsymbol{\phi}_0) - z_t) \right|$$

$$+ \sup_{\|\mathbf{u}\| \leq T} \left| \frac{1}{n} \sum_{t=1}^{n-p} \lambda'(F_z(z_t)) f_z(z_{t,n}^*(\mathbf{u})) \left( \mathbf{u}' \frac{\partial z_t(\boldsymbol{\phi}_0)}{\partial \boldsymbol{\phi}} \right)^2 - 2|\phi_{0r}|^{-1}\tilde{L}\mathbf{u}'\boldsymbol{\Gamma}_p\mathbf{u} \right|$$

$$+ \sup_{\|\mathbf{u}\| \leq T} \left| \frac{1}{2n\sqrt{n}} \sum_{t=1}^{n-p} \lambda'(F_z(z_t)) f_z(z_{t,n}^*(\mathbf{u})) \mathbf{u}' \frac{\partial z_t(\boldsymbol{\phi}_0)}{\partial \boldsymbol{\phi}} \mathbf{u}' \frac{\partial^2 z_t(\boldsymbol{\phi}_{t,n}^*(\mathbf{u}))}{\partial \boldsymbol{\phi} \partial \boldsymbol{\phi}'} \mathbf{u} \right|,$$

where $f_z$ is the density function for $z_1$, $z_{t,n}^*(\mathbf{u})$ is between $z_t$ and $z_t(\boldsymbol{\phi}_0 + n^{-1/2}\mathbf{u})$ and $\boldsymbol{\phi}_{t,n}^*(\mathbf{u})$ is between $\boldsymbol{\phi}_0$ and $\boldsymbol{\phi}_0 + n^{-1/2}\mathbf{u}$. From (A.3), the first term on the right-hand side is $o_p(1)$ and, since there exists a geometrically decaying, nonnegative, real-valued sequence $\{\ddot{\pi}_k\}_{k=-\infty}^{\infty}$ such that

$$\sup_{\|\mathbf{u}\| \leq T} \left| \mathbf{u}' \frac{\partial^2 z_t(\boldsymbol{\phi}_{t,n}^*(\mathbf{u}))}{\partial \boldsymbol{\phi} \partial \boldsymbol{\phi}'} \mathbf{u} \right| \leq \sum_{k=-\infty}^{\infty} \ddot{\pi}_k |z_{t-k}| \qquad \forall t \in \{1, \dots, n-p\}$$

for all $n$ sufficiently large ([7], Section 3.3), the third term is also $o_p(1)$. Using the uniform continuity of $f_z$, the second term equals

$$\sup_{\|\mathbf{u}\| \leq T} \left| \frac{1}{n} \sum_{t=1}^{n-p} \lambda'(F_z(z_t)) f_z(z_t) \left( \mathbf{u}' \frac{\partial z_t(\boldsymbol{\phi}_0)}{\partial \boldsymbol{\phi}} \right)^2 \right.$$
$$\left. - 2\phi_{0r}^{-2}\left( \int_0^1 f_z(F_z^{-1}(s)) \lambda'(s)\, ds \right) \mathbf{u}'\boldsymbol{\Gamma}_p\mathbf{u} \right| + o_p(1),$$

which is $o_p(1)$ by the ergodic theorem. Therefore, (A.10) and, consequently, (A.6) are $o_p(1)$. Similarly, using the uniform continuity of $\lambda'$, it can be shown that (A.7) is $o_p(1)$. $\square$

LEMMA A.4. *For any $T \in (0, \infty)$, as $n \to \infty$,*

$$\sup_{\|\mathbf{u}\|, \|\mathbf{v}\| \leq T} \left| n^{-1} \sum_{t=1}^{n-p} \mathbf{u}' \left[ \lambda \left( \frac{R_t(\boldsymbol{\phi}_0 + n^{-1/2}\mathbf{v})}{n-p+1} \right) - \lambda(F_z(z_t)) \right] \frac{\partial^2 z_t(\boldsymbol{\phi}_0)}{\partial \boldsymbol{\phi} \, \partial \boldsymbol{\phi}'} \mathbf{u} \right| \xrightarrow{P} 0.$$

PROOF. Note that because $F_z$ is a continuous distribution function, it is also uniformly continuous on $\mathbb{R}$. Using (A.11), the Glivenko–Cantelli theorem and the uniform continuity of $F_z$, for any $\varepsilon, \eta > 0$, it can be shown that there exists an integer $m$ such that

$$\mathrm{P} \left( \sup_{\|\mathbf{v}\| \leq T, t \in \{1, \ldots, n-p-m\}} \left| \frac{R_t(\boldsymbol{\phi}_0 + n^{-1/2}\mathbf{v})}{n-p+1} - F_z(z_t) \right| > \eta \right)$$

$$\leq \mathrm{P} \left( \sup_{\|\mathbf{v}\| \leq T, t \in \{1, \ldots, n-p\}} \left| \frac{R_t(\boldsymbol{\phi}_0 + n^{-1/2}\mathbf{v})}{n-p+1} - F_n \left( z_t \left( \boldsymbol{\phi}_0 + \frac{\mathbf{v}}{\sqrt{n}} \right) \right) \right| > \frac{\eta}{3} \right)$$

$$+ \mathrm{P} \left( \sup_{x \in \mathbb{R}} |F_n(x) - F_z(x)| > \frac{\eta}{3} \right)$$

$$+ \mathrm{P} \left( \sup_{\|\mathbf{v}\| \leq T, t \in \{1, \ldots, n-p-m\}} \left| F_z \left( z_t \left( \boldsymbol{\phi}_0 + \frac{\mathbf{v}}{\sqrt{n}} \right) \right) - F_z(z_t) \right| > \frac{\eta}{3} \right)$$

is less than $\varepsilon$ for all $n$ sufficiently large. Hence,

$$\sup_{\|\mathbf{u}\|, \|\mathbf{v}\| \leq T} \left| \frac{1}{n} \sum_{t=1}^{n-p} \mathbf{u}' \left[ \lambda \left( \frac{R_t(\boldsymbol{\phi}_0 + n^{-1/2}\mathbf{v})}{n-p+1} \right) - \lambda(F_z(z_t)) \right] \frac{\partial^2 z_t(\boldsymbol{\phi}_0)}{\partial \boldsymbol{\phi} \, \partial \boldsymbol{\phi}'} \mathbf{u} \right|$$

$$= \sup_{\|\mathbf{u}\|, \|\mathbf{v}\| \leq T} \left| \frac{1}{n} \sum_{t=1}^{n-p} \mathbf{u}' \lambda'(F_{t,n}^*(\mathbf{v})) \frac{\partial^2 z_t(\boldsymbol{\phi}_0)}{\partial \boldsymbol{\phi} \, \partial \boldsymbol{\phi}'} \left[ \frac{R_t(\boldsymbol{\phi}_0 + n^{-1/2}\mathbf{v})}{n-p+1} - F_z(z_t) \right] \mathbf{u} \right|$$

$$\xrightarrow{P} 0,$$

where $F_{t,n}^*(\mathbf{v})$ is between $F_z(z_t)$ and $R_t(\boldsymbol{\phi}_0 + n^{-1/2}\mathbf{v})/(n-p+1)$, since

$$\sup_{\|\mathbf{u}\|, \|\mathbf{v}\| \leq T} \frac{1}{n} \sum_{t=1}^{n-p} \left| \mathbf{u}' \lambda'(F_{t,n}^*(\mathbf{v})) \frac{\partial^2 z_t(\boldsymbol{\phi}_0)}{\partial \boldsymbol{\phi} \, \partial \boldsymbol{\phi}'} \mathbf{u} \right| = O_p(1). \qquad \square$$

For $\mathbf{u} \in \mathbb{R}^p$ and $\delta_1, \delta_2 \in [0, 1]$, let

$$U_n(\mathbf{u}, \delta_1, \delta_2) = \sum_{t=1}^{n-p} \lambda \left( \frac{R_t(\boldsymbol{\phi}_0 + n^{-1/2}\delta_1 \mathbf{u})}{n-p+1} \right) \left[ z_t \left( \boldsymbol{\phi}_0 + \frac{\delta_2 \mathbf{u}}{\sqrt{n}} \right) - z_t \left( \boldsymbol{\phi}_0 + \frac{\delta_1 \mathbf{u}}{\sqrt{n}} \right) \right]$$

and

$$V_n(\mathbf{u}, \delta_1, \delta_2) = \sum_{t=1}^{n-p} \lambda\left(\frac{R_t(\boldsymbol{\phi}_0 + n^{-1/2}\delta_2\mathbf{u})}{n-p+1}\right)\left[z_t\left(\boldsymbol{\phi}_0 + \frac{\delta_2\mathbf{u}}{\sqrt{n}}\right) - z_t\left(\boldsymbol{\phi}_0 + \frac{\delta_1\mathbf{u}}{\sqrt{n}}\right)\right].$$

Using Taylor series expansions,

$$
\begin{aligned}
U_n&(\mathbf{u}, \delta_1, \delta_2)\\
&= \sum_{t=1}^{n-p} \lambda\left(\frac{R_t(\boldsymbol{\phi}_0 + n^{-1/2}\delta_1\mathbf{u})}{n-p+1}\right)\\
&\quad \times \left\{\left[z_t\left(\boldsymbol{\phi}_0 + \frac{\delta_2\mathbf{u}}{\sqrt{n}}\right) - z_t(\boldsymbol{\phi}_0)\right] - \left[z_t\left(\boldsymbol{\phi}_0 + \frac{\delta_1\mathbf{u}}{\sqrt{n}}\right) - z_t(\boldsymbol{\phi}_0)\right]\right\}\\
&= \frac{\delta_2 - \delta_1}{\sqrt{n}}\sum_{t=1}^{n-p}\mathbf{u}'\lambda\left(\frac{R_t(\boldsymbol{\phi}_0 + n^{-1/2}\delta_1\mathbf{u})}{n-p+1}\right)\frac{\partial z_t(\boldsymbol{\phi}_0)}{\partial\boldsymbol{\phi}}\\
&\quad + \frac{1}{2}\frac{\delta_2^2 - \delta_1^2}{n}\sum_{t=1}^{n-p}\mathbf{u}'\lambda\left(\frac{R_t(\boldsymbol{\phi}_0 + n^{-1/2}\delta_1\mathbf{u})}{n-p+1}\right)\frac{\partial^2 z_t(\boldsymbol{\phi}_0)}{\partial\boldsymbol{\phi}\,\partial\boldsymbol{\phi}'}\mathbf{u}\\
&\quad + \frac{1}{2}\frac{\delta_2^2}{n}\sum_{t=1}^{n-p}\mathbf{u}'\lambda\left(\frac{R_t(\boldsymbol{\phi}_0 + n^{-1/2}\delta_1\mathbf{u})}{n-p+1}\right)\\
&\qquad \times \left[\frac{\partial^2 z_t(\boldsymbol{\phi}_n^*(\mathbf{u},\delta_1,\delta_2))}{\partial\boldsymbol{\phi}\,\partial\boldsymbol{\phi}'} - \frac{\partial^2 z_t(\boldsymbol{\phi}_0)}{\partial\boldsymbol{\phi}\,\partial\boldsymbol{\phi}'}\right]\mathbf{u}\\
&\quad - \frac{1}{2}\frac{\delta_1^2}{n}\sum_{t=1}^{n-p}\mathbf{u}'\lambda\left(\frac{R_t(\boldsymbol{\phi}_0 + n^{-1/2}\delta_1\mathbf{u})}{n-p+1}\right)\\
&\qquad \times \left[\frac{\partial^2 z_t(\boldsymbol{\phi}_n^*(\mathbf{u},\delta_1,\delta_1))}{\partial\boldsymbol{\phi}\,\partial\boldsymbol{\phi}'} - \frac{\partial^2 z_t(\boldsymbol{\phi}_0)}{\partial\boldsymbol{\phi}\,\partial\boldsymbol{\phi}'}\right]\mathbf{u}
\end{aligned}
$$

(A.12)

and, similarly,

$$
\begin{aligned}
V_n(\mathbf{u}, \delta_1, \delta_2) &= \frac{\delta_2 - \delta_1}{\sqrt{n}}\sum_{t=1}^{n-p}\mathbf{u}'\lambda\left(\frac{R_t(\boldsymbol{\phi}_0 + n^{-1/2}\delta_2\mathbf{u})}{n-p+1}\right)\frac{\partial z_t(\boldsymbol{\phi}_0)}{\partial\boldsymbol{\phi}}\\
&\quad + \frac{1}{2}\frac{\delta_2^2 - \delta_1^2}{n}\sum_{t=1}^{n-p}\mathbf{u}'\lambda\left(\frac{R_t(\boldsymbol{\phi}_0 + n^{-1/2}\delta_2\mathbf{u})}{n-p+1}\right)\frac{\partial^2 z_t(\boldsymbol{\phi}_0)}{\partial\boldsymbol{\phi}\,\partial\boldsymbol{\phi}'}\mathbf{u}\\
&\quad + \frac{1}{2}\frac{\delta_2^2}{n}\sum_{t=1}^{n-p}\mathbf{u}'\lambda\left(\frac{R_t(\boldsymbol{\phi}_0 + n^{-1/2}\delta_2\mathbf{u})}{n-p+1}\right)\\
&\qquad \times \left[\frac{\partial^2 z_t(\boldsymbol{\phi}_n^*(\mathbf{u},\delta_2,\delta_2))}{\partial\boldsymbol{\phi}\,\partial\boldsymbol{\phi}'} - \frac{\partial^2 z_t(\boldsymbol{\phi}_0)}{\partial\boldsymbol{\phi}\,\partial\boldsymbol{\phi}'}\right]\mathbf{u}
\end{aligned}
$$

(A.13)

$$- \frac{1}{2} \frac{\delta_1^2}{n} \sum_{t=1}^{n-p} \mathbf{u}' \lambda \left( \frac{R_t(\boldsymbol{\phi}_0 + n^{-1/2} \delta_2 \mathbf{u})}{n-p+1} \right)$$

$$\times \left[ \frac{\partial^2 z_t(\boldsymbol{\phi}_n^*(\mathbf{u}, \delta_2, \delta_1))}{\partial \boldsymbol{\phi} \, \partial \boldsymbol{\phi}'} - \frac{\partial^2 z_t(\boldsymbol{\phi}_0)}{\partial \boldsymbol{\phi} \, \partial \boldsymbol{\phi}'} \right] \mathbf{u},$$

where the values of $\boldsymbol{\phi}_n^*(\mathbf{u}, \cdot, \cdot)$ lie between $\boldsymbol{\phi}_0$ and $\boldsymbol{\phi}_0 + n^{-1/2} \mathbf{u}$.

LEMMA A.5.  *For* $\mathbf{u} \in \mathbb{R}^p$, *let* $S_n(\mathbf{u}) = D(\boldsymbol{\phi}_0 + n^{-1/2} \mathbf{u}) - D(\boldsymbol{\phi}_0)$ *and* $S(\mathbf{u}) = \mathbf{u}' \mathbf{N} + |\phi_{0r}|^{-1} (\sigma^2 \tilde{L} - \tilde{K}) \mathbf{u}' \sigma^{-2} \boldsymbol{\Gamma}_p \mathbf{u}$, *where* $\mathbf{N} \sim \mathrm{N}(\mathbf{0}, 2\phi_{0r}^{-2} \{\sigma^2 \tilde{J} - \tilde{K}^2\} \sigma^{-2} \boldsymbol{\Gamma}_p)$. *Then* $S_n(\cdot) \xrightarrow{d} S(\cdot)$ *on* $C(\mathbb{R}^p)$, *the space of continuous functions on* $\mathbb{R}^p$ *where convergence is equivalent to uniform convergence on every compact set.*

PROOF.  Let $\mathbf{u} \in \mathbb{R}^p$ and suppose that $m$ is any positive integer. Because

$$D(\boldsymbol{\phi}_0 + n^{-1/2} \mathbf{u}) - D(\boldsymbol{\phi}_0) = \sum_{k=1}^{m} \left[ D\left( \boldsymbol{\phi}_0 + \frac{k\mathbf{u}}{m\sqrt{n}} \right) - D\left( \boldsymbol{\phi}_0 + \frac{(k-1)\mathbf{u}}{m\sqrt{n}} \right) \right],$$

we have

$$(A.14) \quad \begin{aligned} &\sum_{k=1}^{m} U_n\left( \mathbf{u}, \frac{k-1}{m}, \frac{k}{m} \right) \\ &\leq D(\boldsymbol{\phi}_0 + n^{-1/2} \mathbf{u}) - D(\boldsymbol{\phi}_0) \leq \sum_{k=1}^{m} V_n\left( \mathbf{u}, \frac{k-1}{m}, \frac{k}{m} \right) \end{aligned}$$

by Theorem 2.1. Using (A.12), (A.13) and Lemmas A.1, A.2, A.3 and A.4,

$$\begin{bmatrix} U_n\left( \mathbf{u}, 0, \dfrac{1}{m} \right) \\ U_n\left( \mathbf{u}, \dfrac{1}{m}, \dfrac{2}{m} \right) \\ \vdots \\ U_n\left( \mathbf{u}, \dfrac{m-1}{m}, 1 \right) \\ V_n\left( \mathbf{u}, 0, \dfrac{1}{m} \right) \\ V_n\left( \mathbf{u}, \dfrac{1}{m}, \dfrac{2}{m} \right) \\ \vdots \\ V_n\left( \mathbf{u}, \dfrac{m-1}{m}, 1 \right) \end{bmatrix}$$

$$\overset{d}{\to} \begin{bmatrix} \frac{1}{m}\mathbf{u}'\mathbf{N} - |\phi_{0r}|^{-1}\left[\left(\frac{1}{m}\right)^2 - \left(\frac{0}{m}\right)^2\right]\tilde{K}\mathbf{u}'\sigma^{-2}\boldsymbol{\Gamma}_p\mathbf{u} \\ \frac{1}{m}\mathbf{u}'\mathbf{N} + |\phi_{0r}|^{-1}\left(2\frac{1}{m^2}\sigma^2\tilde{L} - \left[\left(\frac{2}{m}\right)^2 - \left(\frac{1}{m}\right)^2\right]\tilde{K}\right)\mathbf{u}'\sigma^{-2}\boldsymbol{\Gamma}_p\mathbf{u} \\ \vdots \\ \frac{1}{m}\mathbf{u}'\mathbf{N} + |\phi_{0r}|^{-1}\left(2\frac{m-1}{m^2}\sigma^2\tilde{L} - \left[\left(\frac{m}{m}\right)^2 - \left(\frac{m-1}{m}\right)^2\right]\tilde{K}\right)\mathbf{u}'\sigma^{-2}\boldsymbol{\Gamma}_p\mathbf{u} \\ \frac{1}{m}\mathbf{u}'\mathbf{N} + |\phi_{0r}|^{-1}\left(2\frac{1}{m^2}\sigma^2\tilde{L} - \left[\left(\frac{1}{m}\right)^2 - \left(\frac{0}{m}\right)^2\right]\tilde{K}\right)\mathbf{u}'\sigma^{-2}\boldsymbol{\Gamma}_p\mathbf{u} \\ \frac{1}{m}\mathbf{u}'\mathbf{N} + |\phi_{0r}|^{-1}\left(2\frac{2}{m^2}\sigma^2\tilde{L} - \left[\left(\frac{2}{m}\right)^2 - \left(\frac{1}{m}\right)^2\right]\tilde{K}\right)\mathbf{u}'\sigma^{-2}\boldsymbol{\Gamma}_p\mathbf{u} \\ \vdots \\ \frac{1}{m}\mathbf{u}'\mathbf{N} + |\phi_{0r}|^{-1}\left(2\frac{m}{m^2}\sigma^2\tilde{L} - \left[\left(\frac{m}{m}\right)^2 - \left(\frac{m-1}{m}\right)^2\right]\tilde{K}\right)\mathbf{u}'\sigma^{-2}\boldsymbol{\Gamma}_p\mathbf{u} \end{bmatrix}$$

on $\mathbb{R}^{2m}$, since

$$\sup_{\|\mathbf{u}\|,\|\mathbf{v}\|\le T}\frac{1}{n}\sum_{t=1}^{n-p}\left|\mathbf{u}'\left(\frac{\partial^2 z_t(\boldsymbol{\phi}_0 + n^{-1/2}\mathbf{v})}{\partial\boldsymbol{\phi}\,\partial\boldsymbol{\phi}'} - \frac{\partial^2 z_t(\boldsymbol{\phi}_0)}{\partial\boldsymbol{\phi}\,\partial\boldsymbol{\phi}'}\right)\mathbf{u}\right| \overset{P}{\to} 0$$

for any $T > 0$ and $\sup_{s\in(0,1)}|\lambda(s)| < \infty$ by the uniform continuity of $\lambda'$. Hence,

$$\begin{bmatrix} \sum_{k=1}^{m} U_n\left(\mathbf{u}, \frac{k-1}{m}, \frac{k}{m}\right) \\ \sum_{k=1}^{m} V_n\left(\mathbf{u}, \frac{k-1}{m}, \frac{k}{m}\right) \end{bmatrix} \overset{d}{\to} \begin{bmatrix} \mathbf{u}'\mathbf{N} + |\phi_{0r}|^{-1}\left(\frac{m-1}{m}\sigma^2\tilde{L} - \tilde{K}\right)\mathbf{u}'\sigma^{-2}\boldsymbol{\Gamma}_p\mathbf{u} \\ \mathbf{u}'\mathbf{N} + |\phi_{0r}|^{-1}\left(\frac{m+1}{m}\sigma^2\tilde{L} - \tilde{K}\right)\mathbf{u}'\sigma^{-2}\boldsymbol{\Gamma}_p\mathbf{u} \end{bmatrix}$$

on $\mathbb{R}^2$. For any $\varepsilon > 0$, there exists an integer $m$ sufficiently large so that

$$\mathbf{u}'\mathbf{N} + |\phi_{0r}|^{-1}\left(\frac{m-1}{m}\sigma^2\tilde{L} - \tilde{K}\right)\mathbf{u}'\sigma^{-2}\boldsymbol{\Gamma}_p\mathbf{u}$$

and

$$\mathbf{u}'\mathbf{N} + |\phi_{0r}|^{-1}\left(\frac{m+1}{m}\sigma^2\tilde{L} - \tilde{K}\right)\mathbf{u}'\sigma^{-2}\boldsymbol{\Gamma}_p\mathbf{u}$$

are both in an $\varepsilon$-neighborhood of $S(\mathbf{u}) = \mathbf{u}'\mathbf{N} + |\phi_{0r}|^{-1}(\sigma^2\tilde{L} - \tilde{K})\mathbf{u}'\sigma^{-2}\boldsymbol{\Gamma}_p\mathbf{u}$. Thus, for any $\mathbf{u}\in\mathbb{R}^p$, $S_n(\mathbf{u}) \overset{d}{\to} S(\mathbf{u})$. It can be shown similarly that all finite-dimensional distributions of $S_n(\cdot)$ converge to those of $S(\cdot)$.

Also using (A.14), it can be shown that

$$\lim_{\delta\to 0^+}\limsup_{n\to\infty}\mathrm{P}\left(\sup_{\mathbf{u},\mathbf{v}\in K,\|\mathbf{u}-\mathbf{v}\|\le\delta}|S_n(\mathbf{u}) - S_n(\mathbf{v})| > \eta\right) = 0$$

for any $\eta > 0$ and any compact subset $K \subset \mathbb{R}^p$ (see [3], pages 84–86). It follows that $S_n(\cdot)$ must be tight on $C(K)$ and, therefore, because compact $K \subset \mathbb{R}^p$ is arbitrary, $S_n(\cdot) \xrightarrow{d} S(\cdot)$ on $C(\mathbb{R}^p)$ by Theorem 7.1 in [4]. $\quad \square$

Lemma A.6. *If $\varepsilon > 0$ is sufficiently small so that $\boldsymbol{\phi}$ forms a causal polynomial for all $\boldsymbol{\phi} \in \boldsymbol{\Phi} := \{\boldsymbol{\phi} \in \mathbb{R}^p : \|\boldsymbol{\phi} - \boldsymbol{\phi}_0\| \leq \varepsilon\}$, then $n^{-1} \sum_{t=1}^{n-p} z_t^2(\boldsymbol{\phi}) \xrightarrow{a.s.} \mathrm{E}\{\tilde{z}_1^2(\boldsymbol{\phi})\}$ and $n^{-1} D(\boldsymbol{\phi}) \xrightarrow{a.s.} \int_0^1 F_{\tilde{z}(\boldsymbol{\phi})}^{-1}(s)\lambda(s)\,ds$ uniformly on $\boldsymbol{\Phi}$, where $\tilde{z}_t(\boldsymbol{\phi}) := -\phi^{-1}(B^{-1})\phi(B)X_{t+p}$ and $F_{\tilde{z}(\boldsymbol{\phi})}(\cdot)$ is the distribution function for $\tilde{z}_1(\boldsymbol{\phi})$.*

Proof. For any $\boldsymbol{\phi} \in \boldsymbol{\Phi}$, $n^{-1} \sum_{t=1}^{n-p} z_t^2(\boldsymbol{\phi}) \xrightarrow{a.s.} \mathrm{E}\{\tilde{z}_1^2(\boldsymbol{\phi})\}$ and $n^{-1} D(\boldsymbol{\phi}) \xrightarrow{a.s.} \mathrm{E}\{\lambda(F_{\tilde{z}(\boldsymbol{\phi})}(\tilde{z}_1(\boldsymbol{\phi})))\tilde{z}_1(\boldsymbol{\phi})\} = \int_0^1 F_{\tilde{z}(\boldsymbol{\phi})}^{-1}(s)\lambda(s)\,ds$, by the ergodic theorem. Therefore, since $n^{-1} \sum_{t=1}^{n-p} z_t^2(\cdot)$ and $n^{-1} D(\cdot)$ are equicontinuous and uniformly bounded on $\boldsymbol{\Phi}$ almost surely (see Lemma 15 on page 86 of [3]; similar results are obtained in the proof of Proposition 1 in [6]), the lemma follows by the Arzelà–Ascoli theorem. $\quad \square$

**Acknowledgments.** We would like to thank two reviewers and an Associate Editor for their helpful comments. In particular, we are grateful to a reviewer for suggestions that led to Remark 6. We would also like to thank Professor Keh-Shin Lii for supplying the water gun wavelet used in Section 4.2. The work reported here was developed in part under STAR Research Assistance Agreement CR-829095 awarded by the U.S. Environmental Protection Agency (EPA) to Colorado State University. This manuscript has not been formally reviewed by EPA. The views expressed here are solely those of the authors. EPA does not endorse any products or commercial services mentioned in this report.

## REFERENCES

[1] Allal, J., Kaaouachi, A. and Paindaveine, D. (2001). *R*-estimation for ARMA models. *J. Nonparametr. Statist.* **13** 815–831. MR1893753

[2] Andrews, B., Davis, R. A. and Breidt, F. J. (2006). Maximum likelihood estimation for all-pass time series models. *J. Multivariate Anal.* **97** 1638–1659. MR2256234

[3] Andrews, M. E. (2003). Parameter estimation for all-pass time series models. Ph.D. dissertation, Dept. Statistics, Colorado State Univ.

[4] Billingsley, P. (1999). *Convergence of Probability Measures*, 2nd ed. Wiley, New York. MR1700749

[5] Breidt, F. J. and Davis, R. A. (1992). Time-reversibility, identifiability and independence of innovations for stationary time series. *J. Time Ser. Anal.* **13** 377–390. MR1183152

[6] Breidt, F. J., Davis, R. A. and Trindade, A. A. (2001). Least absolute deviation estimation for all-pass time series models. *Ann. Statist.* **29** 919–946. MR1869234

[7] Brockwell, P. J. and Davis, R. A. (1991). *Time Series*: *Theory and Methods*, 2nd ed. Springer, New York. MR1093459

[8] CHI, C.-Y. and KUNG, J.-Y. (1995). A new identification algorithm for allpass systems by higher-order statistics. *Signal Processing* **41** 239–256.

[9] CHIEN, H.-M., YANG, H.-L. and CHI, C.-Y. (1997). Parametric cumulant based phase estimation of 1-D and 2-D nonminimum phase systems by allpass filtering. *IEEE Trans. Signal Processing* **45** 1742–1762.

[10] DAVIS, R. A., KNIGHT, K. and LIU, J. (1992). $M$-estimation for autoregressions with infinite variance. *Stochastic Process. Appl.* **40** 145–180. MR1145464

[11] GIANNAKIS, G. B. and SWAMI, A. (1990). On estimating noncausal nonminimum phase ARMA models of non-Gaussian processes. *IEEE Trans. Acoust. Speech Signal Process.* **38** 478–495. MR1045718

[12] HALLIN, M. (1994). On the Pitman non-admissibility of correlogram-based methods. *J. Time Ser. Anal.* **15** 607–611. MR1312324

[13] HUANG, J. and PAWITAN, Y. (2000). Quasi-likelihood estimation of non-invertible moving average processes. *Scand. J. Statist.* **27** 689–702. MR1804170

[14] JAECKEL, L. A. (1972). Estimating regression coefficients by minimizing the dispersion of the residuals. *Ann. Math. Statist.* **43** 1449–1458. MR0348930

[15] JUREČKOVÁ, J. and SEN, P. K. (1996). *Robust Statistical Procedures*: *Asymptotics and Interrelations*. Wiley, New York. MR1387346

[16] KOUL, H. L. and OSSIANDER, M. (1994). Weak convergence of randomly weighted dependent residual empiricals with applications to autoregression. *Ann. Statist.* **22** 540–562. MR1272098

[17] KOUL, H. L. and SALEH, A. K. MD. E. (1993). $R$-estimation of the parameters of autoregressive [Ar($p$)] models. *Ann. Statist.* **21** 534–551. MR1212192

[18] KOUL, H. L., SIEVERS, G. L. and MCKEAN, J. W. (1987). An estimator of the scale parameter for the rank analysis of linear models under general score functions. *Scand. J. Statist.* **14** 131–141. MR0913258

[19] LII, K.-S. and ROSENBLATT, M. (1988). Nonminimum phase non-Gaussian deconvolution. *J. Multivariate Anal.* **27** 359–374. MR0970960

[20] MUKHERJEE, K. and BAI, Z. D. (2002). $R$-estimation in autoregression with square-integrable score function. *J. Multivariate Anal.* **81** 167–186. MR1901212

[21] ROBINSON, P. M. (1987). Time series residuals with application to probability density estimation. *J. Time Ser. Anal.* **8** 329–344. MR0903762

[22] SILVERMAN, B. W. (1986). *Density Estimation for Statistics and Data Analysis*. Chapman and Hall, London. MR0848134

[23] TERPSTRA, J. T., MCKEAN, J. W. and NARANJO, J. D. (2001). Weighted Wilcoxon estimates for autoregression. *Aust. N. Z. J. Stat.* **43** 399–419. MR1872200

B. ANDREWS
DEPARTMENT OF STATISTICS
NORTHWESTERN UNIVERSITY
2006 SHERIDAN ROAD
EVANSTON, ILLINOIS 60208
USA
E-MAIL: bandrews@northwestern.edu

R. A. DAVIS
F. J. BREIDT
DEPARTMENT OF STATISTICS
COLORADO STATE UNIVERSITY
FT. COLLINS, COLORADO 80523
USA
E-MAIL: rdavis@stat.colostate.edu
         jbreidt@stat.colostate.edu